\documentclass{commat}

\title[On a sum of a multiplicative function]{%
    On a sum of a multiplicative function linked to the divisor function over the set of integers $B$-multiple of $5$
    }

\author{%
    Mihoub Bouderbala
    }

\affiliation{%
    \address{%
    Faculty of Science and Technology, University of Djilali Bounaama.  Department of Mathematics and Informatics, FIMA Laboratory.  Ain Defla, Khemis Miliana (44225), Algeria
        }
    \email{%
    mihoub75bouder@gmail.com
    }
    }

\abstract{%
    Let $d(n)$ and $d^{\ast}(n)$ be the numbers of
divisors and the numbers of unitary divisors of the integer $n\geq1$. In this
paper, we prove that
\[
\underset{n\in\mathcal{B}}{\underset{n\leq x}{\sum}}\frac{d(n)}{d^{\ast}%
(n)}=\frac{16\pi%
{{}^2}%
}{123}\underset{p}{\prod}(1-\frac{1}{2p%
{{}^2}%
}+\frac{1}{2p^{3}})x+\mathcal{O}\left(  x^{\frac{\ln8}{\ln10}+\varepsilon
}\right)  ,~\left(  x\geqslant1,~\varepsilon>0\right)  ,
\]
where $\mathcal{B}$ is the set which contains any integer that is not a
multiple of $5$, but some permutations of its digits is a multiple of $5$.
    }



\VOLUME{31}
\YEAR{2023}
\NUMBER{1}
\firstpage{369}
\DOI{https://doi.org/10.46298/cm.10467}

\begin{paper}

\section{Introduction and main result}

A positive integer is called $\mathcal{A}$-multiple of $5$ if a permutation of
its digits is a multiple of $5$, comprising the identity permutation (for
example $50$, $55$, $505$, $5505$, \ldots). A positive integer is called
$\mathcal{B}$-multiple of $5$ if it is not a multiple of $5$, but some
permutations of its digits is multiples of $5$ (for example $51$, $53$, $107$,
$151$, \ldots). For practical reasons, $\mathcal{A}$ represents the set of
all integers $\mathcal{A}$-multiple of $5$, and $\mathcal{B}$ represents the
set of all integers $\mathcal{B}$-multiple of $5$.

In this paper, we will use deep analytic methods to give an asymptotic formula
to the following sum%
\begin{equation}
\underset{n\in\mathcal{B}}{\underset{n\leq x}{\sum}}\frac{d(n)}{d^{\ast}(n)},
\label{1.}%
\end{equation}
where
\[
d(n)=\underset{d\mid n}{%
{\displaystyle\sum}
}1
\quad \text{ and } \quad
d^{\ast}(n)=\underset{(d,n\mid d)=1}{\underset{d\mid n}{%
{\displaystyle\sum}
}}1,~\ (n\text{ denotes a strictly positive integer}).
\]

In order to estimate the sum $(\ref{1.})$ by noting that the function
$\dfrac{d(n)}{d^{\ast}(n)}~$is multiplicative, we first recall by the
following two concepts of the Riemann zeta function:

We have for all $s\in%
\mathbb{C}
$, such that $\operatorname{Re}(s)>1$,%
\[
\zeta(s)=\underset{p}{\prod}(1-\frac{1}{p^{s}})^{-1},
\]
and
\[
\zeta(s)=1+\frac{1}{s-1}-s\int_{1}^{\infty}\frac{\left\{  t\right\}  }%
{t^{s}+1}dt,
\]
where $\{t\}$ denotes the fractional part of the real $t$.

Recall that according to this last form, the function $\zeta~$extends to a
meromorphic function on $Re(s)>0$, which has a simple pole at $s=1$ with
residue $1$ and no other poles. Moreover, if $a$ is a strictly positive
constant, we have, in the region of the plane defined by the inequalities
$\sigma\geq\dfrac{1}{2},~\sigma\geq1-a/\log\left\vert t\right\vert $, and
$~\sigma\leq2$, the following majoration%
\[
\zeta(\sigma+it)\ll\mathcal{O}\left(  \log\left\vert t\right\vert \right)
,~\text{for }|t|~\text{large enough (see }\cite[p.~54-55]{Key2}
\text{)}.
\]

Secondly, we present the first effective formula of Perron (see
\cite[p.~147]{Key3}): Let 
\[
f(s)=\overset{\infty}%
{\underset{n=1}{%
{\displaystyle\sum}
}}\dfrac{a(n)}{n^{s}}
\]
be the Dirichlet series of finite absolute convergence
abscissa $\sigma_{a}$. Then,\ if $x\geq1,~T\geq1~$and $c>\max(0,\sigma_{a})$,
we have the following asymptotic formula%
\[
\underset{n\leq x}{\sum}a(n)=\frac{1}{2\pi i}\int_{c-iT}^{c+iT}f(s)\frac
{x^{s}}{s}ds+\mathcal{O}\left(  x^{c}\underset{n\geq1}{\sum}\frac{\left\vert
a(n)\right\vert }{n^{c}(1+T\left\vert \ln(x/n)\right\vert }\right)  .
\]

In the following, we will present the main result that has been proven:

\begin{theorem} \label{thm:main}
\textit{For any real }$x\geq1$,\textit{\ we have the following asymptotic
formula}%
\[
\underset{n\in\mathcal{B}}{\underset{n\leq x}{\sum}}\frac{d(n)}{d^{\ast}%
(n)}=\frac{16\pi%
{{}^2}%
}{123}\underset{p}{\prod}(1-\frac{1}{2p%
{{}^2}%
}+\frac{1}{2p^{3}})x+\mathcal{O}\left(  x^{\frac{\ln8}{\ln10}+\varepsilon
}\right)  ,
\]
\textit{where}
\[
\frac{\pi%
{{}^2}%
}{6}\underset{p}{\prod}(1-\frac{1}{2p%
{{}^2}%
}+\frac{1}{2p^{3}})\simeq1.4276565\cdots,\text{ and }\varepsilon>0.
\]

\end{theorem}

The proof of the theorem is based on the following lemmas:

\begin{lemma}
Let $q$ be a prime number or $q=1$. So for any real number $x\geq1$, we have
the following asymptotic formula%
\[
\underset{n\leq x}{\sum}\mathcal{D}(qn)=\frac{2q%
{{}^2}%
-q}{2q%
{{}^2}%
-2q+1}\frac{\pi%
{{}^2}%
}{6}\underset{p}{\prod}(1-\frac{1}{2p%
{{}^2}%
}+\frac{1}{2p^{3}})x+\mathcal{O}(x^{%
\frac12
+\varepsilon}),
\]
\textit{where~}$\mathcal{D}(n)=\dfrac{d(n)}{d^{\ast}(n)}$,\textit{\ and
}$\varepsilon$ denotes a positive real number.
\end{lemma}

\begin{proof}
For a prime number $q$ and a complex number $s$ such that $\operatorname{Re}%
(s)>1$, we put%
\[
~f(s)=\underset{n=1}{\overset{\infty}{%
{\displaystyle\sum}
 }}\dfrac{\mathcal{D}(qn)}{n^{s}},
\]
Then, by the product formula Eulerian $\cite[p.230]{Key1}  $, we
get%
\begin{align*}
f(s) &  =\overset{\infty}{\underset{\alpha=0}{\sum}}\underset{(n_{1}%
,q)=1}{\underset{n_{1}=1}{\overset{\infty}{%
{\displaystyle\sum}
 }}}\frac{\mathcal{D}(q^{\alpha+1}n_{1})}{q^{\alpha s}n_{1}^{s}}\\
&  =\overset{\infty}{\underset{\alpha=0}{\sum}}\frac{\alpha+2}{2q^{\alpha s}%
}\underset{(n_{1},q)=1}{\underset{n_{1}=1}{\overset{\infty}{%
{\displaystyle\sum}
 }}}\dfrac{\mathcal{D}(n_{1})}{n_{1}^{s}}\\
&  =\overset{\infty}{\underset{\alpha=0}{\frac{1}{2}\sum}}\left(  \frac
{\alpha+1}{q^{\alpha s}}+\frac{1}{q^{\alpha s}}\right)  \underset
{(p,q)=1}{\underset{p}{%
{\displaystyle\prod}
 }}\left(  1+\overset{\infty}{\underset{k=1}{\sum}}\frac{\mathcal{D}\left(
p^{k}\right)  }{p^{ks}}\right)  ,
\end{align*}
then%
\begin{align*}
f(s) &  =\frac{1}{2}\left(  \left(  \overset{\infty}{\underset{\alpha=0}{\sum
}}\frac{1}{q^{\alpha s}}\right)  ^{2}+\overset{\infty}{\underset{\alpha
=0}{\sum}}\frac{1}{q^{\alpha s}}\right)  \underset{p}{%
{\displaystyle\prod}
 }\left(  1+\overset{\infty}{\underset{k=1}{\sum}}\frac{\mathcal{D}\left(
p^{k}\right)  }{p^{ks}}\right)  \frac{1}{1+\overset{\infty}{\underset
{k=1}{\sum}}\frac{\mathcal{D}\left(  q^{k}\right)  }{q^{ks}}}\\
&  =\frac{1}{2}\left(  \frac{1}{\left(  1-\frac{1}{q^{s}}\right)  ^{2}}%
+\frac{1}{1-\frac{1}{q^{s}}}\right)  \zeta(s)\zeta(2s)\underset{p}{%
{\displaystyle\prod}
 }\left(  1-\frac{1}{2p^{2s}}+\frac{1}{2p^{3s}}\right)  \frac{2\left(
q^{s}-1\right)  ^{2}}{2q^{2s}-2q^{s}+1}\\
&  =\left(  \frac{2q^{2s}-q^{s}}{2q^{2s}-2q^{s}+1}\right)  \zeta
(s)\zeta(2s)\underset{p}{%
{\displaystyle\prod}
 }\left(  1-\frac{1}{2p^{2s}}+\frac{1}{2p^{3s}}\right)  .
\end{align*}
We notice that the function $f(s)$, is convergent if $\operatorname{Re}%
(s)>\dfrac{1}{2},~$where we recall here that $\zeta(s)=1+\dfrac{1}{s-1}-s%
{\displaystyle\int_{1}^{\infty}}
 \dfrac{\left\{  t\right\}  }{t^{s}+1}dt$. According to Perron's formula, for
all $x\geq1$ and $T\geq1,~$we getting%
\begin{equation}
\underset{n\leq x}{\sum}\mathcal{D}(qn)=\frac{1}{2\pi i}\int_{\frac{3}{2}%
-iT}^{\frac{3}{2}+iT}f(s)\frac{x^{s}}{s}ds+\mathcal{O}(\frac{x^{\frac{3}%
{2}+\varepsilon}}{T}),\label{2}%
\end{equation}
such that $\varepsilon$ is a positive real. 

Now, if we choose a
linear contour integral of $s=\dfrac{3}{2}\pm iT$ to $s=\dfrac{1}{2}\pm iT$,
in this case the function $F(s)=f(s)\dfrac{x^{s}}{s}$, admits a simple pole in
$s=1,~$then%
\[
\frac{1}{2\pi i}\left(  \int_{\frac{1}{2}-iT}^{\frac{3}{2}-iT}+\int_{\frac
{3}{2}-iT}^{\frac{3}{2}+iT}+\int_{\frac{3}{2}+iT}^{\frac{1}{2}+iT}+\int
_{\frac{1}{2}+iT}^{\frac{1}{2}-iT}\right)  f(s)\frac{x^{s}}{s}%
ds=\operatorname{Re}s\left[  f(s)\frac{x^{s}}{s},1\right]  .
\]
Note that $\underset{s\rightarrow1}{\lim}\zeta(s)(s-1)=1$, and we can get
immediately%
\[
\operatorname{Re}s\left[  f(s)\frac{x^{s}}{s},1\right]  =\left(  \frac
{2q^{2}-q}{2q^{2}-2q+1}\right)  \frac{\pi%
{{}^2}%
 }{6}\underset{p}{%
{\displaystyle\prod}
 }\left(  1-\frac{1}{2p^{2}}+\frac{1}{2p^{3}}\right)  x,
\]
such that
\[
\frac{\pi%
{{}^2}%
 }{6}\underset{p}{%
{\displaystyle\prod}
 }\left(  1-\frac{1}{2p^{2}}+\frac{1}{2p^{3}}\right)  \simeq1.4276565\ldots
\]

 By taking $T=x$, and $f(s)=\zeta(s)R(s)$, where
\[
R(s)=\left(  \frac{2q^{2s}-q^{s}}{2q^{2s}-2q^{s}+1}\right)  \zeta
(2s)\underset{p}{%
{\displaystyle\prod}
 }\left(  1-\frac{1}{2p^{2s}}+\frac{1}{2p^{3s}}\right)  ,
\]
we obtain
\begin{align*}
&  \left\vert \frac{1}{2\pi i}\left(  \int_{\frac{1}{2}-iT}^{\frac{3}{2}%
-iT}+\int_{\frac{3}{2}+iT}^{\frac{1}{2}+iT}\right)  \zeta(s)R(s)\frac{x^{s}%
}{s}ds\right\vert \\
&  \ll\int_{\frac{1}{2}}^{\frac{3}{2}}\left\vert \zeta(\sigma+iT)R(s)\frac
{x^{\frac{3}{2}}}{T}\right\vert d\sigma\\
&  \ll\frac{x^{\frac{3}{2}+\varepsilon}}{T}=x^{\frac{1}{2}+\varepsilon},
\end{align*}
and
\[
\left\vert \frac{1}{2\pi i}\int_{\frac{1}{2}+iT}^{\frac{1}{2}-iT}%
\zeta(s)R(s)\frac{x^{s}}{s}ds\right\vert \ll\int_{0}^{T}\left\vert \zeta
(\frac{1}{2}+it)R(s)\frac{x^{\frac{1}{2}}}{t}\right\vert dt\ll x^{\frac{1}%
{2}+\varepsilon}.
\]
So by estimate
\[
\left\vert \frac{1}{2\pi i}\left(  \int_{\frac{1}{2}-iT}^{\frac{3}{2}-iT}%
+\int_{\frac{3}{2}+iT}^{\frac{1}{2}+iT}+\int_{\frac{1}{2}+iT}^{\frac{1}{2}%
-iT}\right)  f(s)\frac{x^{s}}{s}ds\right\vert \ll x^{\frac{1}{2}+\varepsilon
}\text{ },
\]
and from the formula $(\ref{2})$, we get%
\[
\underset{n\leq x}{\sum}\mathcal{D}(qn)=\left(  \frac{2q%
{{}^2}%
 -q}{2q%
{{}^2}%
 -2q+1}\right)  \left(  \frac{\pi%
{{}^2}%
 }{6}\right)  \underset{p}{\prod}(1-\frac{1}{2p%
{{}^2}%
 }+\frac{1}{2p^{3}})x+\mathcal{O}(x^{\frac{1}{2}+\varepsilon}).
\qedhere\]
\end{proof}

\begin{lemma}
For any real $x\geq1$, we have the following asymptotic formula%
\begin{equation}
\underset{n\in\mathcal{A}}{\underset{n\leq x}{\sum}}\mathcal{D}(n)=\frac{\pi%
{{}^2}%
 }{6}\underset{p}{\prod}(1-\frac{1}{2p%
{{}^2}%
}+\frac{1}{2p^{3}})x+\mathcal{O}(x^{\frac{\ln8}{\ln10}+\varepsilon}).
\label{3}%
\end{equation}
\textit{where~}$\mathcal{D}(n)=\dfrac{d(n)}{d^{\ast}(n)}$,\textit{\ and
}$\varepsilon$ denotes a positive real number.
\end{lemma}

\begin{proof}
For any real $x\geq1$, there exists a positive integer $k$ such
that $10^{k}\leq x\leq10^{k+1}$. Consequently, $k\leq\log x\leq k+1$. 
According to the definition of the set $\mathcal{A}$, we know that the number
of integers $\left(  \leq x\right)  $ that is not in $\mathcal{A}$ is
$8^{k+1}$. Indeed, there are $8$ integers composed of a single number, they
are $1,2,3,4,6,7,8,9;$ there are $8%
{{}^2}%
$ integers composed of two digits; and the number of integers composed of
$k$ digits is $8^{k}$. Since%
\[
8^{k}\leq8^{\log x}=x^{\frac{\ln8}{\ln10}},
\]
we get%
\[
\underset{n\notin\mathcal{A}}{\underset{n\leq x}{\sum}}1\leq8+8^{2}%
+8^{3}+...+8^{k+1}\leq\frac{8^{k+2}}{7}\leq\frac{64}{7}8^{k}\leq\frac{64}%
{7}x^{\frac{\ln8}{\ln10}},
\]
Note that for any $\varepsilon>0$ and for all$~n\geq1,~$we have $d(n)\ll
n^{\varepsilon}$, and since $\dfrac{d(n)}{d^{\ast}(n)}\leq d(n)$, we get
$\dfrac{d(n)}{d^{\ast}(n)}\ll n^{\varepsilon}$. Now we apply the lemma
$2~$with $q=1$, we get%
\begin{align*}
\underset{n\in\mathcal{A}}{\underset{n\leq x}{\sum}}\mathcal{D}(n)  &
=\underset{n\leq x}{\sum}\mathcal{D}(n)-\underset{n\notin\mathcal{A}%
}{\underset{n\leq x}{\sum}}\mathcal{D}(n)\\
&  =\underset{n\leq x}{\sum}\mathcal{D}(n)+\mathcal{O}\left(  \underset
{n\notin A}{\underset{n\leq x}{\sum}}x^{\varepsilon}\right) \\
&  =\underset{n\leq x}{\sum}\mathcal{D}(n)+\mathcal{O}\left(  x^{\frac{\ln
8}{\ln10}+\varepsilon}\right) \\
&  =\frac{\pi%
{{}^2}%
}{6}\underset{p}{\prod}(1-\frac{1}{2p%
{{}^2}%
}+\frac{1}{2p^{3}})x+\mathcal{O}\left(  x^{\frac{\ln8}{\ln10}+\varepsilon
}\right)  .
\end{align*}
This proves Lemma $3$.
\end{proof}

\section{Proof of Theorem~\ref{thm:main}}

In this section, we complete the proof of Theorem. From the definition of the
set $\mathcal{A}$ and set $\mathcal{B}$, we know the relation between them.
Therefore%
\begin{align*}
\underset{n\in\mathcal{B}}{\underset{n\leq x}{\sum}}\mathcal{D}(n)  &
=\underset{n\in\mathcal{A}}{\underset{n\leq x}{\sum}}\mathcal{D}%
(n)-\underset{5n\leq x}{\sum}\mathcal{D}(5n)\\
&  =\underset{n\in\mathcal{A}}{\underset{n\leq x}{\sum}}\mathcal{D}%
(n)-\underset{n\leq\frac{x}{5}}{\sum}\mathcal{D}(5n).
\end{align*}
Now we use the two results of Lemmas $2$ and $3$, we get%
\begin{align*}
\underset{n\in\mathcal{B}}{\underset{n\leq x}{\sum}}\frac{d(n)}{d^{\ast}(n)}
&  =\frac{\pi%
{{}^2}%
}{6}\underset{p}{\prod}(1-\frac{1}{2p%
{{}^2}%
}+\frac{1}{2p^{3}})x-\frac{3\pi^{2}}{82}\underset{p}{\prod}(1-\frac{1}{2p%
{{}^2}%
}+\frac{1}{2p^{3}})x+\mathcal{O}\left(  x^{\frac{\ln8}{\ln10}+\varepsilon
}\right) \\
&  =\frac{16\pi%
{{}^2}%
}{123}\underset{p}{\prod}(1-\frac{1}{2p%
{{}^2}%
}+\frac{1}{2p^{3}})x+\mathcal{O}\left(  x^{\frac{\ln8}{\ln10}+\varepsilon
}\right)  .
\end{align*}
This completes the proof of the theorem.

\subsection*{Acknowledgments}
The author thank Professor Olivier Bordell\`{e}s for his help and attention to
this work.


\EditInfo{July 23, 2021}{September 06, 2022}{Attila Bérczes}

\end{paper}